\documentclass{article}

\usepackage[
  lang = american, 
 paper = hardcover, openaccess,
]{ems-ecm}

\numberwithin{equation}{section}

\newtheorem{theorem}{Theorem}[section]

\newcommand{\rt}{\tau_{\vartheta}}
\def\Blue#1{#1}
\def\Red#1{#1}

\def\aone{a}
\def\bone{b}
\def\betil{\tilde{\beta}}
\renewcommand{\Re}{\textrm{Re}}
\newcommand{\uld}[1]{\underline{d#1}}

\definecolor{grey}{rgb}{0.5, 0.5, 0.5}
\usepackage{comment,tikz,pgfplots,bm}

\numberwithin{equation}{section}

\begin{document}


\title{Forward and Inverse Problems in Nonlinear Acoustics}
\titlemark{Forward and Inverse Problems in Nonlinear Acoustics}



\emsauthor{1}{
	\givenname{Barbara}
	\surname{Kaltenbacher}
	\mrid{616341} %
	\orcid{0000-0002-3295-6977}}{B.~Kaltenbacher} 
	
\Emsaffil{1}{
	\department{Department of Mathematics}
	\organisation{University of Klagenfurt}
	\rorid{01a2bcd34} 
	\address{Universit\"atsstra\ss e 65-67}
	\zip{9020}
	\city{Klagenfurt}
	\country{Austria}
	\affemail{barbara.kaltenbacher@aau.at}}
%
%
	
\classification[YYyYY]{XXxXX}


\begin{abstract}
The importance of ultrasound is well established in the imaging of human tissue. In order to enhance image quality by exploiting nonlinear effects, recently techniques such as harmonic imaging and nonlinearity parameter tomography have been put forward. 
As soon as the pressure amplitude exceeds a certain bound, the classical linear wave equation loses its validity and more general nonlinear versions have to be used.
Another characteristic property of ultrasound propagation in human tissue is frequency power law attenuation, leading to fractional derivative damping models in time domain.
In this contribution we will first of all dwell on modeling nonlinearity on the one hand and fractional damping on the other hand. 
Moreover we will give an idea on the challenges in the analysis of the resulting PDEs and discuss some parameter asymptotics.
Finally, we address a relevant inverse problems in this context, the above mentioned task of nonlinearity parameter imaging, which leads to a coefficient identification problem for a quasilinear wave equation.

\end{abstract}

\maketitle


\section{Introduction}

\begin{figure}
\includegraphics[width=0.78\textwidth]{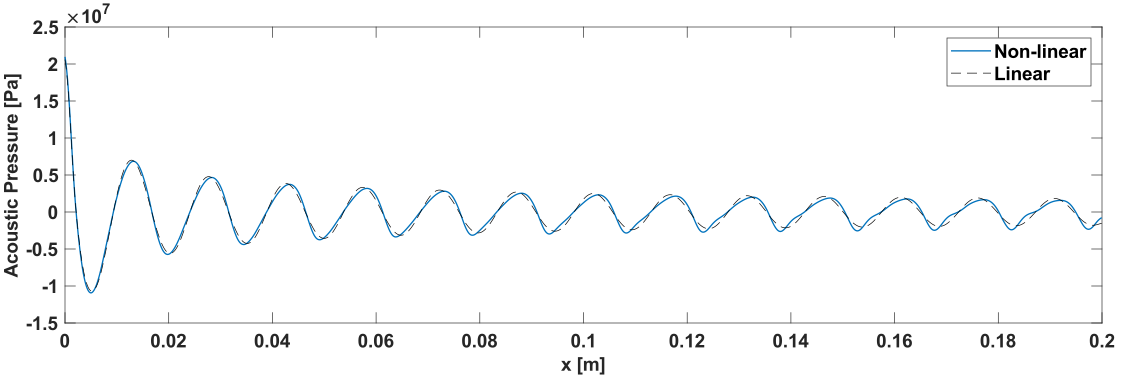}
\includegraphics[width=0.58\textwidth]{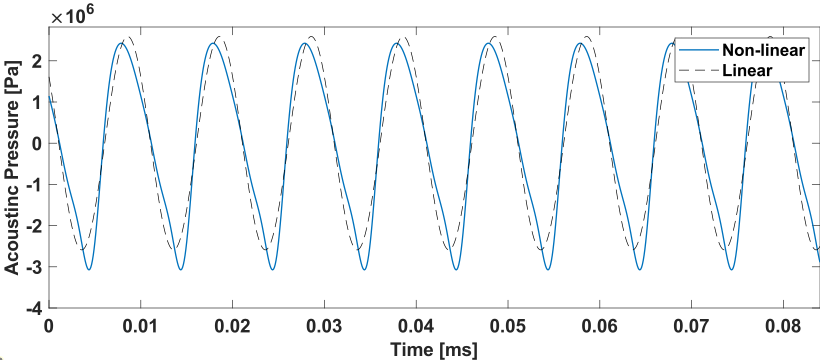}
\caption{Linear versus nonlinear wave propagation in space (top) and time (bottom) domain.
\label{fig:sawteeth}
}
\end{figure}

\begin{figure}
\begin{center}
\includegraphics[width=0.58\textwidth]{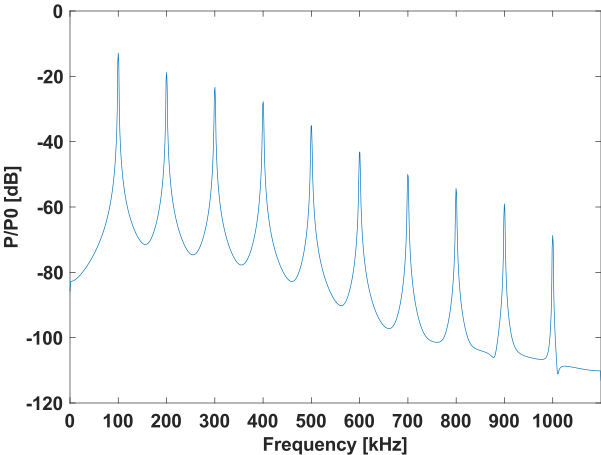}
\end{center}
\caption{Higher harmonics in frequency domain.
\label{fig:harmonics}
}
\end{figure}

The difference between linear and nonlinear acoustic wave propagation becomes apparent by a steepening of waves, see Figure~\ref{fig:sawteeth}. This is due to a pressure dependent wave speed, which is higher in compressed areas and lower in low density regions. This fact is visible in the wave profiles from both a space and from a time perspective, and can be also be read off from the PDEs, cf. \eqref{wavespeed} below. Another way of characterizing nonlinearity is by looking at contributions in frequency domain in case of a sinusiodal excitation, which in the linear case results in a response at the excitation frequency only, whereas in the nonlinear case, so-called higher harmonics are generated, see Figure~\ref{fig:harmonics}. Also this effect is very well visible in the governing mathematical models, cf. \eqref{multiharmonics}.

{}
\section{Models of nonlinear acoustics}\label{sec:models}
\subsection{Physical Principles}
The main physical quantities in acoustics are\\[1ex]
\begin{tabular}{ll}
$\bullet$ acoustic particle velocity {$\mathbf{v}$};    & $\bullet$ absolute temperature ${\vartheta}$;\\
$\bullet$ acoustic pressure {$p$};  &  $\bullet$ heat flux ${\boldsymbol{q}}$;
\\
$\bullet$ mass density {$\varrho$};  & $\bullet$ entropy ${\eta}$;
\end{tabular}
\\[1ex]
These are usually decomposed into mean (zero order) and fluctuating (first order) parts
\[
\mathbf{v}=\mathbf{v}_0+\mathbf{v}_\sim=\mathbf{v}\,, \qquad p=p_0+p_\sim\,, \qquad \varrho=\varrho_0+\varrho_\sim\,, \qquad\text{etc.}
\]
and subject to the following physical laws
\begin{itemize}
\item momentum conservation = 
Navier Stokes 
equation {(with $\nabla\times\mathbf{v}=0$)}: 
\[
{\varrho} \Big({\mathbf{v}_t} + \nabla({\mathbf{v}} \cdot {\mathbf{v}})\Big) + \nabla {p} =
\Big(\tfrac{4\mu_V}{3} + \zeta_V \Big) \Delta {\mathbf{v}}
\]
\item mass conservation = equation of continuity: 
\qquad ${\displaystyle 
{\varrho_t}+\nabla\cdot({\varrho}{\mathbf{v}})=0
}$
\item entropy equation: 
\qquad ${\displaystyle 
{\varrho}{\vartheta}({\eta}_t+{\mathbf{v}}\cdot\nabla{\eta})=-\nabla\cdot{\boldsymbol{q}}
}$
\item equation of state:
\qquad ${\displaystyle 
{\varrho_\sim} = \frac{1}{c^2}{p_\sim} - \frac{1}{\varrho_0 c^4}\tfrac{B}{2A} {p_\sim}^2
- \frac{a}{\varrho_0 c^4}\Big( \frac{1}{c_V} - \frac{1}{c_p} \Big)  {{p_\sim}_t}
}$
\item Gibbs equation: 
\qquad ${\displaystyle 
{\vartheta} d{\eta} =c_V d{\vartheta}-{p}\frac{1}{\varrho^2}d{\varrho}
}$
\end{itemize}
where 
$c_p$ / $c_V$ is the specific heat at constant pressure / volume, and 
$\zeta_V$ / $\mu_V$ the bulk / shear viscosity.
The missing 6th equation in this system governing the 6 unknown functions ${\mathbf{v}}$, ${p}$, ${\varrho}$, ${\vartheta}$, ${\boldsymbol{q}}$, ${\eta}$ is the relation between temperature and heat flux, which (like the pressure -- density relation) is subject to constitutive modeling.
Classically, Fourier's law 
\begin{equation}\label{Fourier}
    {\boldsymbol{q}}=-K\nabla{\vartheta}
\end{equation}
(with {$K$ denoting the thermal conductivity}) has been employed for this purpose. However this is known to lead to the infinite speed of propagation paradox.
For this reason, a relaxation term has been introduced, that leads to the 
Maxwell-Cattaneo law 
\begin{equation}\label{MaxwellCattaneo}
\tau {\boldsymbol{q}}_t+{\boldsymbol{q}}=-K\nabla{\vartheta}, 
\end{equation}
which contains the relaxation time {$\tau$} as an additional parameter and
allows for so-called ``thermal waves'' (also called the second sound phenomenon).

\subsection{Classical Models of Nonlinear Acoustics}\label{sec:ClassMod}
In the linear case, where quadratic and higher order terms (e.g., products of first order terms) are skipped, it suffices to consider momentum and mass conservation along with the equation of state, and the velocity can be eliminated by taking the divergence in momentum conservation and subtracting the time derivative of mass conservation. Upon substitution of the (linearized) pressure -- density relation, this leads to the well-known second order wave equation $p_{tt}-c^2\Delta p=0$. 
In the nonlinear case, the same procedure, but retaining quadratic terms, 
(that is, application of the so-called Blackstock's scheme in nonlinear acoustics) leads to the classical models of nonlinear acoustics.
\begin{itemize}
\item {\bf Kuznetsov's equation} \cite{LesserSeebass68,Kuznetsov71}
\[
{
{p_\sim}_{tt}- c^2\Delta p_\sim - \delta \Delta {p_\sim}_t =
\left(\tfrac{B}{2A \varrho_0 c^2} p_\sim^2 + \varrho_0 |\mathbf{v}|^2\right)_{tt}
}
\]
where 
$\varrho_0\mathbf{v}_t=-\nabla p$ 
for the 
particle velocity $\mathbf{v}$
and the 
pressure $p$, i.e., 
\begin{equation}\label{Kuznetsov}
{
\psi_{tt}-c^2\Delta\psi-\delta\Delta\psi_t=\Bigl(\tfrac{B}{2A\, c^2}(\psi_t)^2
+|\nabla\psi|^2\Bigr)_t,
}
\end{equation}
since $\nabla\times\mathbf{v}=0$ hence 
{
$\mathbf{v}=-\nabla\psi$
} 
for some 
velocity potential $\psi$
\item The {\bf Westervelt equation} \cite{Westervelt63}, which is obtained from Kuznetsov's equation via the approximation $\varrho_0|\mathbf{v}|^2\approx\frac{1}{\varrho_0c^2}({p_\sim})^2$ that corresponds to neglecting non-cumulative nonlinear effects 
\begin{equation}\label{Westervelt}
{
{p_\sim}_{tt}- c^2\Delta p_\sim - \delta \Delta {p_\sim}_t =
\frac{1}{\varrho_0 c^2}\Bigl(1+\tfrac{B}{2A}\Bigr) {p_\sim^2}_{tt}
}
\end{equation}
\end{itemize}
Here $\delta=h (\text{Pr}(\tfrac43+\tfrac{\zeta_V}{\mu_V})+\gamma-1)$ is the diffusivity of sound; (containing the thermal diffusivity $h$ and the Prandtl number Pr) and 
$\tfrac{B}{A}$ 
is the  nonlinearity parameter. 

\subsection{Advanced Models of Nonlinear Acoustics}
Taking into account further physical effects requires to consider temperature, heat flux and entropy as variables as well.
Modeling in nonlinear acoustics is a highly active field and only a small portion of the resulting PDE models has been subject to rigorous mathematical investigation so far.
We here mention a few examples of these.
\begin{itemize}
\item The {\bf  Jordan-Moore-Gibson-Thompson} JMGT equation \cite{Christov2009,JordanMaxwellCattaneo14,Straughan2010}
\begin{equation}\label{JMGT}
\tau\psi_{ttt}+\psi_{tt}-c^2\Delta \psi - (\delta+\tau c^2)\Delta \psi_t
= \left(\tfrac{B}{2A c^2}(\psi_t)^2+|\nabla\psi|^2\right)_t
\end{equation}
where $\tau$ is the relaxation time, results from replacing Fourier's law \eqref{Fourier} of heat conduction by \eqref{MaxwellCattaneo}; later on, we will also encounter fractional versions thereof.
The mentioned second sound phenomenon can be seen by considering the auxiliary quantity 
$z:= \psi_t + \frac{c^2}{\delta+\tau c^2}\psi$, which solves a weakly damped wave equation
\[
z_{tt}  - \tilde{c} \Delta z +\gamma z_t
=r(z,\psi)
\]
with $\tilde{c}=c^2+\tfrac{\delta}{\tau}$, $\gamma=\tfrac{1}{\tau}- \tfrac{c^2}{\delta+\tau c^2}>0$ and a lower order term $r(z,\psi)$.
\item The {\bf Blackstock-Crighton} BCBJ equation \cite{Blackstock63,BrunnhuberJordan15,Crighton79}
\begin{equation}\label{BCBJ}
\left(\partial_t -a\Delta\right)
\left(\psi_{tt}-c^2\Delta\psi - \delta\Delta\psi_t \right)-ra\Delta\psi_t
=\left(\tfrac{B}{2Ac^2}(\psi_t^2)+\left|\nabla\psi\right|^2\right)_{tt}
\end{equation}
where $a=\frac{\nu}{\mbox{Pr}}$ is the thermal diffusivity and $r$ another constant.
\end{itemize}
Both models can be related to the classical Kuznetsov equation \eqref{Kuznetsov} by -- first of all just formally -- setting the parameters $\tau$ and $a$ to zero.
We will take a look at the mathematical justification of these parameter limits in Section~\ref{sec:limits}.

The amount of literature on the analysis of these and some other advanced models of nonlinear acoustics is vast and we do not attempt to provide an overview that could claim completeness in any sense. 
Some selected references are \cite{bucci2018feedback,DekkersRozanova,dell2016moore,Fritz,KL12_Kuznetsov,KLP12_JordanMooreGibson,KT18_ModelsNlAcoustics,Mizohata1993global,Nikolic15}.

Just to get an idea of the mathematical challenges arising in the analysis (referring to the above quoted literature for details) let us take a look at the question of well-posedness of a prototypical initial-boundary value problem 
for the Westervelt equation as the most simple of the above mentioned models 
\begin{align*}
u_{tt}- c^2\Delta u - b \Delta u_t &=
\tfrac{\kappa}{2} (u^2)_{tt}&& \mbox{ in }\Omega\\
\frac{\partial u}{\partial n}&=g&&\mbox{ on }\partial\Omega\\
u(t=0)=u_0\,, \ u_t(t=0)&=u_1&&\mbox{ in }\Omega,
\end{align*}
where we use the mathematically common notation $u$ for the PDE solution, which plays the physical role of the acoustic pressure.
Differentiating out the quadratic term and a slight rearrangement leads to
\[
{(1-\kappa u)}u_{tt}- c^2\Delta u - b \Delta u_t =
\kappa(u_t)^2 
\]
which reveals the fact that the coefficient of the second time derivative depends on $u$ itself.
This leads to potential {degeneracy for $u\geq \frac{1}{\kappa}$} and similarly affects the other models (Kuznetsov, Jordan-Moore-Gibson-Thompson, Blackstock-Crighton).
The typical approach of proving well-posedness therefore relies on smallness of the state and proceeds as follows:
\begin{itemize} 
\item  employ energy estimates to obtain a bound on $u$ in $C(0,T;H^2(\Omega))$;
\item  use smallness of $u$ in $C(0,T;H^2(\Omega))$ and the embedding $H^2(\Omega)\to L^\infty(\Omega)$ to guarantee 
$1-\kappa u\geq\underline{\alpha}>0$;
\item  apply a combination of these arguments in a fixed point scheme.
\end{itemize}
This also illustrates {state dependence of the effective wave speed}, since a further rearrangement allows to rewrite the equation as 
\[
u_{tt}-\tilde{c}^2\Delta  u - \tilde{b}(u)\Delta  u_t
=f(u)
\]
with 
\begin{equation}\label{wavespeed}
\tilde{c}(u)=\frac{c}{\sqrt{1-\kappa u}}, \quad \tilde{b}(u)=\frac{b}{1-\kappa u},\quad f(u)=\frac{\kappa(u_t)^2}{1-\kappa u},
\end{equation}
as long as $1-\kappa u>0$ (otherwise the equation loses its validity as a wave propagation model).

\section{Parameter asymptotics}\label{sec:limits}
In this section we point to a few analytical results on limits that interrelate the models in Section~\ref{sec:models}. We will be explicit about function space settings and mathematical tools only in part for the first example (the limit as $\tau\searrow$ in the JMGT equation, cf. Section~\ref{sec:limittau}). The other limits will be discussed briefly and collectively in Section~\ref{sec:limitsother}; details can be found in the cited papers.

\subsection{Vanishing relaxation time in the 
Jordan-Moore-Gibson-Thompson equation} \label{sec:limittau}
Considering a family of solutions  $(\psi^\tau)_{\tau\in(0,\overline{\tau}]}$ to 
\[
\tau\psi^\tau_{ttt}+\psi^\tau_{tt}-c^2\Delta \psi^\tau - b\Delta \psi^\tau_t
= \left(\tfrac{B}{2A c^2}(\psi^\tau_t)^2+|\nabla\psi^\tau|^2\right)_t
\]
with $b=\delta+\tau c^2$, the typical question to be answered is, whether and in which function spaces a limit $\psi^0$ of $\psi^\tau$ as $\tau\searrow0$ exists and under which conditions it can be shown to solve Kuznetsov's equation \eqref{Kuznetsov}.
\\[1ex]
Some comments are in order.
\begin{itemize}
\item
We will consider the ``Westervelt type'' and the ``Kuznetsov type'' equation; without and with the gradient nonlinearity $|\nabla \psi|^2_t$.
\item
For $\tau=0$ (classical Westervelt and Kuznetsov equation) the reformulation of the linearization as a first order system leads to an analytic semigroup and maximal parabolic regularity.
These properties do not hold any more with $\tau>0$, that is, the equation loses its ``parabolic nature''.
This is consistent with physics, corresponding to a transition from infinite to finite propagation speed. As $\tau\to0$ the PDE changes from hyperbolic to parabolic. 
\item 
As in the classical models \eqref{Kuznetsov}, \eqref{Westervelt}, potential degeneracy can be an issue, since we have 
\[
\begin{aligned}
\tau\psi^\tau_{ttt}+\psi^\tau_{tt}-c^2\Delta \psi^\tau - b\Delta \psi^\tau_t
=& \left(\frac{\kappa}{2}(\psi^\tau_t)^2+|\nabla\psi^\tau|^2\right)_t\\
=& \kappa\psi^\tau_t\psi^\tau_{tt}+|\nabla\psi^\tau|^2_t\\
\Longleftrightarrow\quad
\tau\psi^\tau_{ttt}+{(1-\kappa\psi^\tau_t)}\psi^\tau_{tt}-c^2\Delta \psi^\tau &- b\Delta \psi^\tau_t
=|\nabla\psi^\tau|^2_t.
\end{aligned}
\]
\end{itemize}
The plan of the analysis is to first establish well-posedness of the linearized equation along with energy estimates and use these results to prove well-posedness of the Westervelt or Kuznetsov type JMGT equation for $\tau>0$ by a fixed point argument.
Relying on uniform in $\tau$ bounds obtained in the energy analysis and (weak) compactness of bounded sets in the underlying solution spaces, we then take limits as $\tau\to0$. 
We refer to \cite{Bongarti2021,JMGT,JMGT_Neumann} for details and only sketch the above steps in a concrete function space setting for the Westervelt type case, while mentioning in passing that an extension to the Kuznetsov type JMGT equation requires higher order energy estimates.

We start by considering the linearized problem with variable coefficient $\alpha$
\begin{equation} \label{ibvp_linear}
\begin{aligned}
\begin{cases}
{\tau}\psi_{ttt}+\alpha(x,t)\psi_{tt}-c^2\Delta \psi - b\Delta \psi_t = f \quad \mbox{ in }\Omega\times(0,T), \\[1mm]
\psi=0 \quad \mbox{ on } \partial \Omega\times(0,T),\\[1mm]
(\psi, \psi_t, \psi_{tt})=(\psi_0, \psi_1, \psi_2) \quad \mbox{ in }\Omega\times \{0\}.
\end{cases}
\end{aligned}
\end{equation}
Under the assumptions 
\begin{align} \label{non-degeneracy_assumption}
\alpha(x,t)\geq\underline{\alpha}>0\ \mbox{ on }\Omega \ \ \mbox{  a.e. in } \Omega \times (0,T).
\end{align}
\begin{equation}\label{eq:alphagammaf_reg_Wes}
\begin{aligned} 
&\alpha \in 
L^\infty(0,T; L^\infty(\Omega))\cap L^\infty(0,T; W^{1,3}(\Omega)), \\
&f\in H^1(0,T; L^2(\Omega)).
\end{aligned}
\end{equation}
\begin{equation}\label{reg_init}
(\psi_0, \psi_1, \psi_2)\in 
H_0^1(\Omega) \cap H^2(\Omega)\times H_0^1(\Omega) \cap H^2(\Omega)\times H_0^1(\Omega).
\end{equation} 
we obtain the following result
\begin{theorem}\label{thm:lin}
	Let $c^2$, $b$, $\tau>0$, and let $T>0$. Let the assumptions \eqref{non-degeneracy_assumption}, \eqref{eq:alphagammaf_reg_Wes}, \eqref{reg_init} hold. Then there exists a unique  solution
	\begin{equation*}
	\psi \in  X^W:=W^{1, \infty}(0,T;H_0^1(\Omega) \cap H^2(\Omega)) \cap W^{2, \infty}(0,T; H_0^1(\Omega))
	\cap H^3(0,T; L^2(\Omega)).
	\end{equation*}
The solution fullfils the estimate 
	\begin{equation*}
	\begin{aligned}
\|\psi\|_{W,\tau}^2:=
&  {\tau}^2 \|\psi_{ttt}\|^2_{L^2L^2} + {\tau} \|\psi_{tt}\|^2_{L^\infty H^1}+\|\psi_{tt}\|^2_{L^2 H^1}+\|\psi\|^2_{W^{1,\infty} H^2}
\\[1mm]
\leq&\,C(\alpha, T, {\tau})\left(|\psi_0|^2_{H^2}+|\psi_1|^2_{H^2}+{\tau}|\psi_2|^2_{H^1}+ \| f\|^2_{L^\infty L^2}+\|f_t\|^2_{L^2 L^2}\right). 
	\end{aligned}
	\end{equation*}
If additionally 
$\|\nabla \alpha\|_{L^\infty L^3}<\frac{\underline{\alpha}}{C^\Omega_{H^1,L^6} }$
holds, then $C(\alpha, T, {\tau})$ is independent of ${\tau}$.
\end{theorem}
\textit{Idea of proof:} \\ 
As standard for proving well-posedness of evolutionary PDEs, we combine a Galerkin discretization with energy estimates and weak limits; however, due to the need for higher regularity, 
we derive energy estimates by nonstandard testing; in this case we test by $-\Delta\psi_{tt}$ and with $\tau\psi_{ttt}$.
\hfill 
    $\diamondsuit$

\medskip

This allows to prove well-posedness of the Westervelt type JMGT equation
	\begin{equation*} \label{ibvp_Westervelt_MC}
	\begin{aligned}
	\begin{cases}
	{\tau}\psi_{ttt}+(1-\kappa\psi_t)\psi_{tt}-c^2\Delta \psi - b\Delta \psi_t = 0\,
 \quad \mbox{ in }\Omega\times(0,T), \\[1mm]
	\psi=0 \quad \mbox{ on } \partial \Omega\times(0,T),\\[1mm]
	(\psi, \psi_t, \psi_{tt})=(\psi_0, \psi_1, \psi_2) \quad \mbox{ in }\Omega\times \{0\},
	\end{cases}
	\end{aligned}
	\end{equation*}
\begin{theorem} \label{th:wellposedness_Wes}
	Let $c^2$, $b>0$, $\kappa\in\mathbb{R}$ and $T>0$. There exist $\rho$,$\rho_0>0$ such that for all $(\psi_0,\psi_1,\psi_2)\in H_0^1(\Omega) \cap H^2(\Omega)\times H_0^1(\Omega) \cap H^2(\Omega)\times H_0^1(\Omega)
$ satisfying 
\begin{equation*}
\|\psi_0\|_{H^2(\Omega)}^2+\|\psi_1\|_{H^2(\Omega)}^2+{\tau}\|\psi_2\|_{H^1(\Omega)}^2\leq\rho_0^2\,,
\end{equation*} 
there exists a unique solution $\psi\in X^W$ and $\|\psi\|_{W,\tau}^2\leq\rho^2$.
\end{theorem}
\textit{Idea of proof:} \\ 
We apply Banach's Contraction Principle to the fixed point operator $\mathcal{T}:\phi\mapsto\psi$ solving \eqref{ibvp_linear} with $\alpha=1-\kappa\phi_t$, $f=0$.
Invariance of $\mathcal{T}$ on $B_\rho^{X^W}$ can be obtained from the energy estimate in Theorem~\ref{thm:lin}. 
Also contractivity $\|\mathcal{T}(\phi_1)-\mathcal{T}(\phi_2)\|_{W,\tau}\leq q \|\phi_1-\phi_2\|_{W,\tau}$ can be established by means of the estimate in Theorem~\ref{thm:lin}, since 
$\hat{\psi}=\psi_1-\psi_2=\mathcal{T}(\phi_1)-\mathcal{T}(\phi_2)$ solves \eqref{ibvp_linear} with 
$\alpha=1-\kappa\phi_{1\,t}$ and $f=\kappa\hat{\phi}_t \psi_{2\,tt}$,
where $\hat{\phi}=\phi_1-\phi_2$.
\hfill 
    $\diamondsuit$

\medskip

In order to take limits for vanishing relaxation time, we 
consider the $\tau$-independent part of the norms
\[	\begin{aligned}
&\|\psi\|_{W,\tau}^2:=  {\tau}^2 \|\psi_{ttt}\|^2_{L^2L^2} + {\tau} \|\psi_{tt}\|^2_{L^\infty H^1}+\|\psi_{tt}\|^2_{L^2 H^1}+\|\psi\|^2_{W^{1,\infty} H^2},
\end{aligned}
\]
that is,
\begin{equation*}
\begin{aligned}
&\|\psi\|_{\bar{X}^W}^2:=\|\psi_{tt}\|^2_{L^2 H^1}+\|\psi\|^2_{W^{1,\infty} H^2}\,, 
\end{aligned}
\end{equation*}
since these norms are uniformly bounded, independently of $\tau$,
and impose smallness of initial data in the space  
\[
\begin{aligned}
&X_0^W:=H_0^1(\Omega)\cap H^2(\Omega)\times H_0^1(\Omega)\cap H^2(\Omega) \times H_0^1(\Omega).
\end{aligned}
\]
This allows us to prove the following, cf~\cite{JMGT}.
\begin{theorem} \label{th:limits}
	Let $c^2$, $b$, $\kappa \in \mathbb{R}$ and $T>0$. Then there exist $\bar{\tau}$, $\rho_0>0$ such that for all $(\psi_0,\psi_1,\psi_2)\in B_{\rho_0}^{X_0^W}$, the family $(\psi^\tau)_{\tau\in(0,\bar{\tau})}$ of solutions to the Westervelt type JMGT equation converges weakly* in $\bar{X}^W$ to a solution $\bar{\psi}\in \bar{X}^W$ of the Westervelt equation with initial conditions $\bar{\psi}(0)=\psi_0$, $\bar{\psi}_t(0)=\psi_1$ as $\tau\searrow0$.
\end{theorem}

\begin{figure}
\begin{center}
		\input{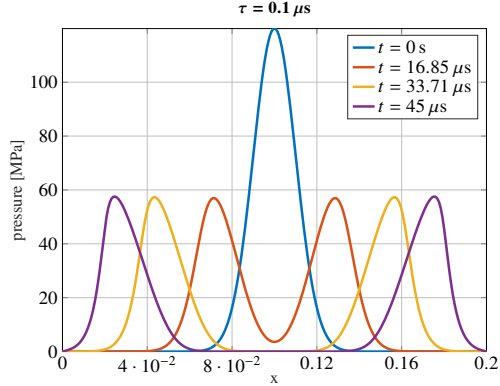}
	\end{center}
\caption{Snapshots of pressure $p=\varrho \psi_t$ for fixed relaxation time $\tau=0.1 \, \mu$s
\label{fig:plot_JMGT_t}}
\end{figure}

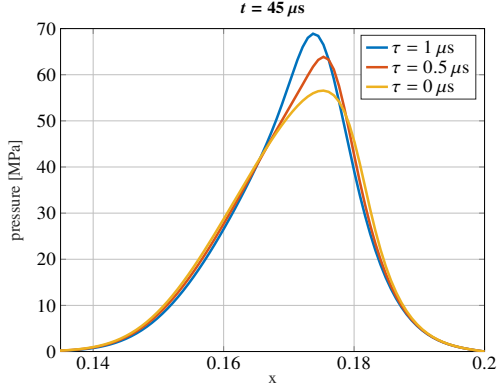
\begin{figure}
	\begin{center}
%
\definecolor{mycolor1}{rgb}{0.00000,0.44700,0.74100}%
\definecolor{mycolor2}{rgb}{0.85000,0.32500,0.09800}%
\definecolor{mycolor3}{rgb}{0.92900,0.69400,0.12500}%
\begin{tikzpicture}[scale=0.475, font=\LARGE]

\begin{axis}[%
width=4.602in,
height=3.566in,
at={(0.772in,0.481in)},
scale only axis,
xmin=0.135,
xmax=0.2, xtick={0.14, 0.16, 0.18, 0.2},
xlabel style={font=\Large\color{white!15!black}},
xlabel={x},
ymin=0,
ymax=70,
ylabel style={font=\Large\color{white!15!black}},
ylabel={pressure [MPa]},
axis background/.style={fill=white},
xmajorgrids,
ymajorgrids,
title style={font=\bfseries},
title={\Large \boldsymbol{$t=45\,\mu$}s},
legend style={legend cell align=left, align=left, draw=white!15!black}
]
\addplot [color=mycolor1, line width=2.0pt]
  table[row sep=crcr]{%
0.134661354581667	0.146586064296613\\
0.135458167330682	0.18835560284225\\
0.136254980079684	0.242349627411315\\
0.137051792828686	0.311432472414978\\
0.137848605577688	0.398920787109162\\
0.138645418326689	0.508593717210971\\
0.139442231075691	0.644685934062437\\
0.140239043824707	0.811860494810134\\
0.141035856573708	1.01515953797782\\
0.14183266932271	1.25993223775913\\
0.142629482071712	1.5517411474158\\
0.143426294820713	1.89624988315489\\
0.144223107569715	2.29909681792516\\
0.145019920318731	2.76576084556142\\
0.145816733067733	3.30142614789561\\
0.146613545816734	3.91085312927024\\
0.147410358565736	4.59826224459182\\
0.148207171314738	5.36723640637396\\
0.149003984063739	6.22064616190237\\
0.149800796812755	7.1606000837083\\
0.150597609561757	8.18842102956577\\
0.151394422310759	9.30464729782888\\
0.15219123505976	10.509056382824\\
0.152988047808762	11.8007081256593\\
0.153784860557764	13.178003616227\\
0.15458167330678	14.6387562687096\\
0.155378486055781	16.1802721140734\\
0.156175298804783	17.7994376324978\\
0.156972111553785	19.4928155954419\\
0.157768924302786	21.2567527581558\\
0.158565737051788	23.0875083621872\\
0.159362549800804	24.9814198966656\\
0.160159362549805	26.9351329062923\\
0.160956175298807	28.9459345675065\\
0.161752988047809	31.0122441218985\\
0.16254980079681	33.1343209903226\\
0.163346613545812	35.3152405671667\\
0.164143426294814	37.5621350141045\\
0.16494023904383	39.8875655838513\\
0.165737051792831	42.3106369998091\\
0.166533864541833	44.8570395641551\\
0.167330677290835	47.5566110479219\\
0.168127490039836	50.4363778078231\\
0.168924302788838	53.5067681848704\\
0.171314741035857	63.2204387767926\\
0.172111553784859	65.9933232286102\\
0.172908366533861	68.0015766726021\\
0.173705179282862	68.893155860443\\
0.174501992031878	68.419937878696\\
0.17529880478088	66.5158529323397\\
0.176095617529882	63.3138543458718\\
0.176892430278883	59.0969952321079\\
0.177689243027885	54.2154532795212\\
0.180079681274904	38.6808501011064\\
0.180876494023906	33.9012765279211\\
0.181673306772907	29.5009416164829\\
0.182470119521909	25.5142506990459\\
0.183266932270911	21.9463266546376\\
0.184063745019927	18.7836996493853\\
0.184860557768928	16.0019971336055\\
0.18565737051793	13.571177047278\\
0.186454183266932	11.4589351493213\\
0.187250996015933	9.6328330770075\\
0.188047808764935	8.06156499757888\\
0.188844621513951	6.71566118165086\\
0.189641434262953	5.56783251743769\\
0.190438247011954	4.59309137879279\\
0.191235059760956	3.76873664786928\\
0.192031872509958	3.07425860348658\\
0.192828685258959	2.49119821906523\\
0.193625498007975	2.00298168558747\\
0.194422310756977	1.59474222641809\\
0.195219123505979	1.2531358172231\\
0.19601593625498	0.966154128970913\\
0.196812749003982	0.722936112308545\\
0.197609561752984	0.513578639756005\\
0.198406374501985	0.328946116757521\\
0.199203187251001	0.160478915635551\\
0.200000000000003	0\\
};
\addlegendentry{$\tau=1 \,\mu$s}

\addplot [color=mycolor2, line width=2.0pt]
  table[row sep=crcr]{%
0.134661354581674	0.172856720498245\\
0.135458167330675	0.223355167723525\\
0.136254980079684	0.287689754065852\\
0.137051792828686	0.3689636538288\\
0.137848605577688	0.470758301143789\\
0.138645418326696	0.597142930632728\\
0.139442231075698	0.752664603556205\\
0.1402390438247	0.942314824253259\\
0.141035856573708	1.17146998119484\\
0.14183266932271	1.44580456399468\\
0.142629482071712	1.77117829076823\\
0.143426294820721	2.15350071562419\\
0.144223107569722	2.59857924885982\\
0.145019920318724	3.11195848009847\\
0.145816733067733	3.69875991456601\\
0.146613545816734	4.36353152591224\\
0.147410358565736	5.11011581080859\\
0.148207171314738	5.94154343588261\\
0.149003984063746	6.85995731119965\\
0.149800796812748	7.8665693690054\\
0.15059760956175	8.96164978156812\\
0.151394422310759	10.1445461438775\\
0.15219123505976	11.4137284488875\\
0.152988047808762	12.7668546078803\\
0.153784860557771	14.200850802321\\
0.154581673306772	15.712000998861\\
0.155378486055774	17.2960404342034\\
0.156175298804783	18.9482485529741\\
0.156972111553785	20.6635377712202\\
0.157768924302786	22.4365352119243\\
0.158565737051795	24.2616554460784\\
0.160159362549798	28.04522222366\\
0.161752988047809	31.967381431434\\
0.164143426294821	38.0063718831429\\
0.167330677290835	46.0924321370293\\
0.168924302788845	50.043451084846\\
0.172111553784859	57.8657083391596\\
0.172908366533868	59.8573535719618\\
0.173705179282869	61.7338451492673\\
0.174501992031871	63.2167883809034\\
0.17529880478088	63.8796671548211\\
0.176095617529882	63.2752748463226\\
0.176892430278883	61.1444505149992\\
0.177689243027892	57.558060013679\\
0.178486055776894	52.8776575157648\\
0.179282868525895	47.5863886139004\\
0.180079681274897	42.1282604912735\\
0.180876494023906	36.8291237350917\\
0.181673306772907	31.8886035480154\\
0.182470119521909	27.4062720372048\\
0.183266932270918	23.4143285651814\\
0.18406374501992	19.9045049511637\\
0.184860557768921	16.8467698261965\\
0.18565737051793	14.2010649562003\\
0.186454183266932	11.9241766938873\\
0.187250996015933	9.97351418933079\\
0.188047808764942	8.309018738396\\
0.188844621513944	6.89400263921642\\
0.189641434262946	5.69539155775256\\
0.190438247011954	4.6836518673806\\
0.191235059760956	3.8325633285635\\
0.192031872509958	3.11892468016689\\
0.192828685258966	2.52223996131303\\
0.193625498007968	2.02440872494373\\
0.19442231075697	1.6094305549778\\
0.195219123505979	1.26312707144896\\
0.19601593625498	0.972880823870206\\
0.196812749003982	0.727388872493371\\
0.197609561752991	0.516428161756579\\
0.198406374501992	0.330629727637863\\
0.199203187250994	0.161259111491646\\
0.200000000000003	0\\
};
\addlegendentry{$\tau=0.5 \, \mu$s}

\addplot [color=mycolor3, line width=2.0pt]
  table[row sep=crcr]{%
0.134661354581674	0.185274706211494\\
0.135458167330675	0.241671214320974\\
0.136254980079684	0.313097443423885\\
0.137051792828686	0.402852550046418\\
0.137848605577688	0.514740283529378\\
0.138645418326696	0.65308097207491\\
0.139442231075698	0.822702938296324\\
0.1402390438247	1.02890842920625\\
0.141035856573708	1.27741009244865\\
0.14183266932271	1.57423568968168\\
0.142629482071712	1.92560119242938\\
0.143426294820721	2.33775534839277\\
0.144223107569722	2.81680204915383\\
0.145019920318724	3.36850971614942\\
0.145816733067733	3.99811915658533\\
0.146613545816734	4.71016222353825\\
0.147410358565736	5.50830320326822\\
0.148207171314738	6.39521281544188\\
0.149003984063746	7.3724817485633\\
0.149800796812748	8.44057680714878\\
0.15059760956175	9.59883910413592\\
0.151394422310759	10.8455202717373\\
0.15219123505976	12.177850342005\\
0.152988047808762	13.5921293358257\\
0.153784860557771	15.0838342674804\\
0.154581673306772	16.647733500972\\
0.155378486055774	18.278001409515\\
0.156175298804783	19.9683275129115\\
0.156972111553785	21.7120155677189\\
0.158565737051795	25.3312639868525\\
0.160159362549798	29.0773331739284\\
0.164940239043823	40.4575127838447\\
0.165737051792831	42.2884887113909\\
0.166533864541833	44.0765281318376\\
0.167330677290835	45.8112142604975\\
0.168127490039844	47.4811192200558\\
0.168924302788845	49.073487874543\\
0.169721115537847	50.5738250983595\\
0.170517928286856	51.9653540110852\\
0.171314741035857	53.2282974721371\\
0.172111553784859	54.338916464448\\
0.172908366533868	55.2682154606617\\
0.173705179282869	55.9801913815294\\
0.174501992031871	56.4294829245964\\
0.17529880478088	56.5582814483026\\
0.176095617529882	56.2925289542742\\
0.176892430278883	55.5380034863114\\
0.177689243027892	54.1786741197141\\
0.178486055776894	52.0840737655347\\
0.179282868525895	49.139662291263\\
0.180079681274897	45.3144103509371\\
0.180876494023906	40.747270826429\\
0.182470119521909	30.8028948126736\\
0.183266932270918	26.171242643418\\
0.18406374501992	22.0427412307828\\
0.184860557768921	18.459076314383\\
0.18565737051793	15.3953896932982\\
0.186454183266932	12.7982119175145\\
0.187250996015933	10.6074432427969\\
0.188047808764942	8.76552185382271\\
0.188844621513944	7.22063255568701\\
0.189641434262946	5.92771222786189\\
0.190438247011954	4.84793576453214\\
0.191235059760956	3.94806181674937\\
0.192031872509958	3.19966000721233\\
0.192828685258966	2.57834591379026\\
0.193625498007968	2.06316901848672\\
0.19442231075697	1.6360397681828\\
0.195219123505979	1.28125985747023\\
0.19601593625498	0.985112934947686\\
0.196812749003982	0.735500955173634\\
0.197609561752991	0.521627740863323\\
0.198406374501992	0.333705538079421\\
0.199203187250994	0.162685480109594\\
0.200000000000003	0\\
};
\addlegendentry{$\tau=0 \,\mu$s}

\end{axis}
\end{tikzpicture}%
	\end{center}
\caption{Pressure wave for different relaxation parameters $\tau$ at final time $t=45\mu$s.
\label{fig:plot_JMGT_tau}}
\end{figure}

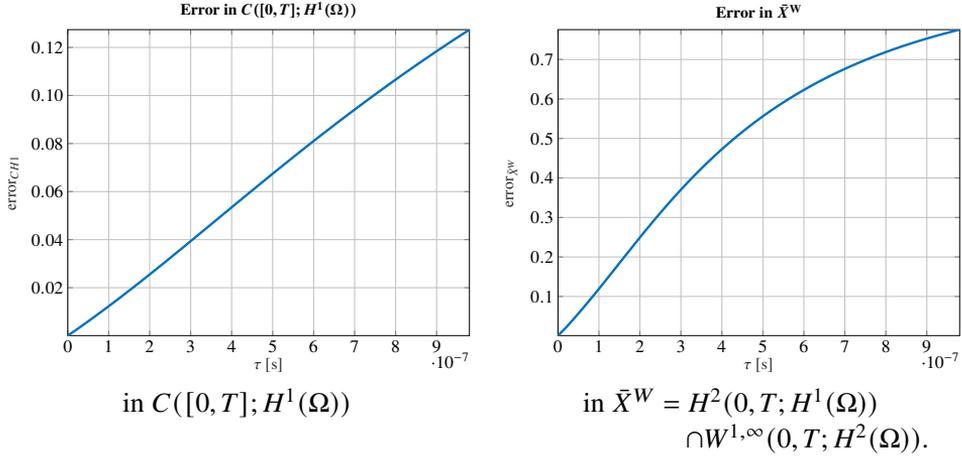
\begin{figure}
\begin{tabular}{cc}
%
\definecolor{mycolor1}{rgb}{0.00000,0.44700,0.74100}%
\begin{tikzpicture}[scale=0.45, font=\LARGE]

\begin{axis}[%
width=4.602in,
height=3.566in,
at={(0.772in,0.481in)},
scale only axis,
xmin=1.00000008274037e-10,
xmax=9.80099999992046e-07,
xlabel style={font=\Large\color{white!15!black}},
xlabel={\Large $\tau \, [\textup{s}]$},
ymin=6.42731803469632e-05,
ymax=0.127400486593906,
ylabel style={at={(-0.03,0.5)}, font=\Large\color{white!15!black}},
ylabel={$\textup{error}_{CH1}$},
axis background/.style={fill=white},
y tick label style={
	/pgf/number format/.cd,
	fixed,
	fixed zerofill,
	precision=2,
	/tikz/.cd
},
xmajorgrids,
ymajorgrids,
title style={font=\bfseries},
title={\Large Error in \boldsymbol{$C([0,T]; H^1(\Omega))$}},
legend style={legend cell align=left, align=left, draw=white!15!black}
]
\addplot [color=mycolor1, line width=2.0pt]
  table[row sep=crcr]{%
9.80099999992046e-07	0.127400486593906\\
9.401000000131e-07	0.12296819671296\\
9.00100000006399e-07	0.118429670668816\\
8.60099999999697e-07	0.113784468687131\\
8.20099999992996e-07	0.109032646955549\\
7.80099999986295e-07	0.104174862535778\\
7.40100000007349e-07	0.0992124962822906\\
7.00100000000647e-07	0.094147796131937\\
6.60099999993946e-07	0.0889840429539007\\
6.20099999987245e-07	0.0837257406294813\\
5.80100000008299e-07	0.0783788308945823\\
5.40100000001598e-07	0.0729509313431745\\
5.00099999994896e-07	0.0674515912854545\\
4.40099999998722e-07	0.0590949804980019\\
3.60100000013075e-07	0.0478333985433927\\
2.80099999999672e-07	0.0365809209297603\\
2.40099999992971e-07	0.0310175442428655\\
2.0009999998627e-07	0.0255306384491289\\
1.60100000007324e-07	0.0201455685299709\\
1.20100000000622e-07	0.014884905004411\\
1.0010000001115e-07	0.0123070054582203\\
8.00999999939211e-08	0.0097664017274092\\
6.01000000044483e-08	0.00726435135496456\\
4.00999999872198e-08	0.00480174410407444\\
2.00999999977469e-08	0.00237927898201032\\
1.00000008274037e-10	6.42731803469632e-05\\
};
\end{axis}
\end{tikzpicture}%
& 
%
\definecolor{mycolor1}{rgb}{0.00000,0.44700,0.74100}%
\begin{tikzpicture}[scale=0.45, font=\LARGE]

\begin{axis}[%
width=4.602in,
height=3.566in,
at={(0.772in,0.481in)},
scale only axis,
xmin=1.00000008274037e-10,
xmax=9.80100000047557e-07,
xlabel style={font=\Large\color{white!15!black}},
xlabel={\Large $\tau \, [\textup{s}]$},
ymin=0.000850673042971595,
ymax=0.775506875805468,
ylabel style={at={(-0.01,0.5)}, font=\Large\color{white!15!black}},
ylabel={$\textup{error}_{\bar{X}^W}$},
axis background/.style={fill=white},
xmajorgrids,
ymajorgrids,
title style={font=\bfseries},
title={\Large Error in \boldsymbol{$\bar{X}^{\textup{W}}$}},
legend style={legend cell align=left, align=left, draw=white!15!black}
]
\addplot [color=mycolor1, line width=2.0pt]
  table[row sep=crcr]{%
9.80100000047557e-07	0.775506875805468\\
9.60099999947062e-07	0.770235030351333\\
9.40099999957589e-07	0.764734330988659\\
9.20099999968116e-07	0.758993809393763\\
9.00099999978643e-07	0.753001914304792\\
8.8009999998917e-07	0.746746481742432\\
8.60099999999697e-07	0.740214704568574\\
8.40100000010224e-07	0.733393101622989\\
8.20100000020751e-07	0.726267486737512\\
8.00100000031279e-07	0.71882293802669\\
7.80100000041806e-07	0.711043767962051\\
7.60100000052333e-07	0.702913494879316\\
7.40099999951838e-07	0.694414816738766\\
7.20099999962365e-07	0.685529588172628\\
7.00099999972892e-07	0.676238802111098\\
6.80099999983419e-07	0.666522577603831\\
6.60099999993946e-07	0.656360155833154\\
6.40100000004473e-07	0.645729906808815\\
6.20100000015e-07	0.634609349799713\\
6.00100000025527e-07	0.622975191288366\\
5.80100000036055e-07	0.610803385073181\\
5.60100000046582e-07	0.598069220210226\\
5.40099999946086e-07	0.584747443727482\\
5.20099999956614e-07	0.570812426556258\\
5.00099999967141e-07	0.556238382915484\\
4.80099999977668e-07	0.540999655459236\\
4.60099999988195e-07	0.525071081047394\\
4.40099999998722e-07	0.508428454643195\\
4.20100000009249e-07	0.491049111824449\\
4.00100000019776e-07	0.472912653698107\\
3.80100000030303e-07	0.454001841376463\\
3.60100000040831e-07	0.434303689049905\\
3.40100000051358e-07	0.41381078617108\\
3.20099999950862e-07	0.392522877496462\\
3.0009999996139e-07	0.370448723297579\\
2.80099999971917e-07	0.347608247100277\\
2.60099999982444e-07	0.324034951097278\\
2.40099999992971e-07	0.299778532403848\\
2.20100000003498e-07	0.274907558638198\\
2.00100000014025e-07	0.249511950035212\\
1.60100000035079e-07	0.197623352866078\\
1.20099999945111e-07	0.145293964180653\\
1.00099999955638e-07	0.119412364826836\\
8.00999999661656e-08	0.0939694120943644\\
6.00999999766927e-08	0.0691374911532111\\
4.00999999872198e-08	0.0450620089195959\\
2.00999999977469e-08	0.0218555442874305\\
1.00000008274037e-10	0.000850673042971595\\
};

\end{axis}
\end{tikzpicture}%
\\
		in $C([0,T]; H^1(\Omega))$ & in  $\bar{X}^{W}=H^2(0,T;H^1(\Omega))$\\
								   &\hspace*{2cm}$\cap W^{1,\infty}(0,T;H^2(\Omega))$.
\end{tabular}
\caption{Relative errors as $\tau\to0$
\label{fig:plot_error_CH1}}
\end{figure}

Figures~\ref{fig:plot_JMGT_t}, \ref{fig:plot_JMGT_tau}, and \ref{fig:plot_error_CH1} illustrate these convergence results by showing some numerical simulations comparing Westervelt-JMGT and Westervelt solutions in the following setting.
\begin{itemize}
\item water in a 1-d channel geometry
$$c=1500 \, \textup{m}/\textup{s}, \ \delta=6\cdot 10^{-9} \, \textup{m}^2/\textup{s}, \ \rho = 1000 \, \textup{kg}/\textup{m}^3, \ B/A=5;$$
\item space discretization with B-splines (Isogeometric Analysis): quadratic basis functions, globally $C^2$; 251 dofs on $\Omega=[0,0.2m]$
\item time discretization by Newmark scheme, adapted to 3rd order equation; 800 time steps on $[0,T]=[0,45\mu s]$
\item initial conditions \quad $\left(\psi_0, \psi_1, \psi_2 \right)= \left (0, \, \mathcal{A}\, \text{exp}\left (-\frac{(x-0.1)^2}{2\sigma^2} \right),\, 0 \right)$
with $\mathcal{A}=8\cdot 10^{4} \, \textup{m}^2/\textup{s}^2$. 
\end{itemize}

\subsection{Some further parameter limits}\label{sec:limitsother}
\begin{itemize}
    \item Vanishing diffusivity of sound $\delta\searrow0$ in 
{Kuznetsov's equation \eqref{Kuznetsov} }:
\[
\psi^\delta_{tt}-c^2\Delta\psi^\delta - {\delta}\Delta\psi^\delta_t 
=\left(\tfrac{B}{2Ac^2}(\psi^\delta_t)^2+\left|\nabla\psi^\delta\right|^2\right)_{t}
\]
with an {undamped quasilinear wave equation}
\begin{equation}\label{undampedKuznetsov}
\psi_{tt}-c^2\Delta\psi  
=\left(\tfrac{B}{2Ac^2}(\psi_t)^2+\left|\nabla\psi\right|^2\right)_{t}
\end{equation}
as a limiting problem, cf.~\cite{b2zero}.
The challenge here lies in the fact that $\delta>0$ is crucial for global in time well-posedness and exponential decay  in $d\in\{2,3\}$ space dimensions.
As a byproduct for $\delta=0$ we recover results (in particular on required regularity of initial data) from \cite{DoerflerGernerSchnaubelt2016} on the undamped Westervelt equation and obtain a new result on local in time well-posedness of the undamped Kuznetsov equation \eqref{undampedKuznetsov}. 
A similar limiting analyisis for $\delta\searrow0$  has been carried out for the JMGT equation \eqref{JMGT} in \cite{b2zeroJMGT}.
\item {Vanishing thermal diffusivity in the Blackstock-Crighton equation \eqref{BCBJ}}:\\ 
Here the aim is to recover solutions to Kuznetsov's equation \eqref{Kuznetsov} as a limit for $a\searrow0$ of solutions $\psi^a$ to
\[
\left(\partial_t -{a}\Delta\right)
\left(\psi^a_{tt}-c^2\Delta\psi^a - \delta\Delta\psi^a_t \right)-ra\Delta\psi^a_t
=\left(\tfrac{B}{2Ac^2}({\psi^a_t}^2)+\left|\nabla\psi^a\right|^2\right)_{tt}
\]
cf.~\cite{KT18_ModelsNlAcoustics}.
Note that this process also requires one integration with respect to time. To obtain sufficiently regular solutions, (as needed for the underlying arguments for nondegeneracy),  consistency of the initial data is thus needed
\[\psi_2-c^2\Delta\psi_0 - \delta\Delta\psi_1 
=\tfrac{B}{Ac^2}\psi_1\psi_2+2\nabla\psi_0\cdot\nabla\psi_1.
\]
\item {Limit $\alpha\nearrow 1$ of the differentiation order in fractional JMGT equations}:\\
Finally, we also consider families of time-fractional versions of the Jordan-Moore-Gibson-Thompson equation \eqref{JMGT}, such as 
\[
\tau^{{\alpha}} \partial_t^{2+{\alpha}}\psi^\alpha+\psi^\alpha_{tt}-c^2\Delta \psi^\alpha - (\delta+\tau^{{\alpha}} c^2)\Delta \partial_t^{{\alpha}}\psi^\alpha
=\left(\tfrac{B}{2A c^2}(\psi^\alpha_t)^2+|\nabla\psi^\alpha|^2\right)_t
\]
and convergence of their solutions to one of \eqref{JMGT} (with fixed positive $\tau$), cf.~\cite{fracJMGT}; see also \eqref{fracJMGTeqs} below.
Like in the $\tau\searrow0$ limit, the leading order time derivative changes; however here this happens in a continuous manner. 
One of the key tasks in this work was the derivation of proper models from physical balance and constitutive laws. This brings us to the topic of the next section.
\end{itemize}

\section{Fractional attenuation models in ultrasonics}
\begin{figure}
\begin{center}
\includegraphics[width=0.5\textwidth]{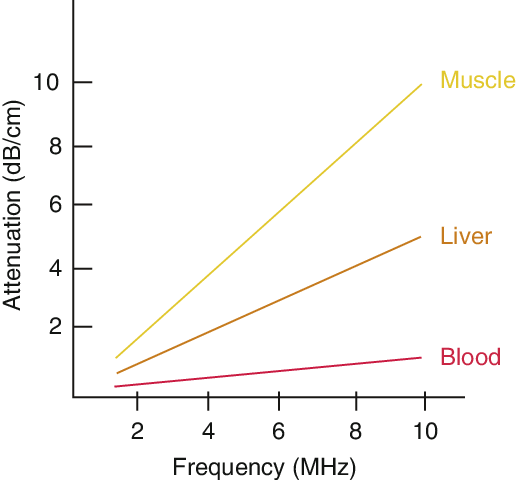}\\
    \caption{Figure 2.6 in [Chan\&Perlas, Basics of Ultrasound Imaging, 2011]
    \label{fig:attUS}}
\end{center}
\end{figure}

As Figure~\ref{fig:attUS} shows, attenuation of ultrasound depends both on the tissue type (which by the way also provides a way of imaging as mentioned in Section~\ref{sec:inverse}) and on the frequency. In view of the semilogarithmic scaling in the graph (the attenuation is given in dezibel per centimeter), the slopes of the lines correspond to  exponents of the frequency in the resulting attenuation law.
Since physical balance equations by their nature only contain integer derivatives, fractional (or more generally,  nonlocal) modeling happens on the level of constitutive laws.
In acoustics this can be done in terms of 
\begin{itemize}
\item pressure -- density relation and/or
\item temperature -- heat flux relation 
\end{itemize}
Below we will partly describe the derivation of such models in the simpler linear setting; nonlinear versions can be obtained by inserting the constitutive equations into the balance laws and applying the procedure mentioned at the beginning of Section~\ref{sec:ClassMod}, leading to second (or higher) order in time wave equations.

As pressure -- density relation often rely on analogy to viscoelasticity (which was one or the historical origins of fractional modeling \cite{Caputo:1967}), we start with a brief recap on fractional modeling there.

\subsection{Fractional models of (linear) viscoelasticity}
In continuum mechanics, the deformation of a body is described by 
\begin{itemize}
\item
the equation of motion (resulting from balance of forces)
\[
\varrho\mathbf{u}_{tt}=\mbox{div}\sigma+\mathbf{f};
\]
\item 
the characterization of strain as the symmetric gradient of displacements (in a geometrically linearized setting)
\[
\epsilon=\frac12(\nabla\mathbf{u}+(\nabla\mathbf{u})^T);
\]
\item
a constitutive model, that is, a stress-strain relation.
\end{itemize}
In here, $\mathbf{u}$ is the vector of displacements,
$\sigma$ and $\epsilon$ are the stress and strain tensors and $\varrho$ is the mass density.

Classical examples of constitutive models are
\[
\begin{aligned}
\mbox{Hooke's law (linear elasticity):}\quad&
\sigma = b_0 \epsilon
\\
\mbox{Newton model:}\quad&
\sigma = b_1 \epsilon_t
\\
\mbox{Kelvin-Voigt model:}\quad&
\sigma = b_0 \epsilon + b_1 \epsilon_t
\\
\mbox{Maxwell model:}\quad&
\sigma + a_1\sigma_t = b_0 \epsilon
\\
\mbox{Zener model:}\quad&
\sigma + a_1\sigma_t = b_0 \epsilon + b_1 \epsilon_t
\end{aligned}
\]
with coefficients $a_i$, $b_i$.
Replacing the first time derivatives in there by fractional ones $\partial_t^{\beta}$ (cf. Section~\ref{sec:fracderiv}) with $\beta\in[0,1)$, one arrives at (time-)fractional versions
\[
\begin{aligned}
\mbox{\Blue{fractional} Newton model:}\quad&
\sigma = b_1 \Blue{\partial_t^{\beta}} \epsilon
\\
\mbox{\Blue{fractional} Kelvin-Voigt model:}\quad&
\sigma = b_0 \epsilon + b_1 \Blue{\partial_t^{\beta}} \epsilon
\\
\mbox{\Blue{fractional} Maxwell model:}\quad&
\sigma + a_1\Blue{\partial_t^{\alpha}}\sigma = b_0 \epsilon
\\
\mbox{\Blue{fractional} Zener model:}\quad&
\sigma + a_1\Blue{\partial_t^{\alpha}}\sigma = b_0 \epsilon + b_1 \Blue{\partial_t^{\beta}}\epsilon
\\
\mbox{which are comprised in the general model class:}\quad&
\sum_{n=0}^N a_n \Blue{\partial_t^{\alpha_n}} \sigma = \sum_{m=0}^M b_m \Blue{\partial_t^{\beta_m}} \epsilon
\end{aligned}
\]
cf. \cite{Atanackovic:2014,Caputo:1967}.

\subsection{Fractional models of (linear) acoustics via pressure -- density relation}
As we have already seen in Section~\ref{sec:models}, the fundamental laws of acoustics are first of all
\begin{itemize}
\item
balance of momentum 
\[
\varrho_0 \mathbf{v}_t = -\nabla p +\mathbf{f}
\]
\item 
balance of mass  
\[
\varrho \nabla\cdot \mathbf{v}=-\varrho_t 
\]
\item combined with an equation of state; 
the most simple linear one is 
\[
\frac{\varrho_\sim}{\varrho_0} = \frac{p_\sim}{p_0} 
\]
and analogously to the general viscoelastic law above, can be generalized to 
\[
\sum_{m=0}^M b_m \Blue{\partial_t^{\beta_m}} \frac{\varrho_\sim}{\varrho_0}
= \sum_{n=0}^N a_n \Blue{\partial_t^{\alpha_n}} \frac{p_\sim}{p_0}. 
\]
\end{itemize}
Inserting the equation of state into a combination of the balance laws that eliminates the velocity leads to fractional acoustic wave equations.
Some of the most commonly used instances are as follows \cite{Holm2019waves,Szabo:1994}.
\begin{itemize}
\item 
Caputo-Wismer-Kelvin wave equation (fractional Kelvin-Voigt):
\begin{equation*}
p_{tt} - b_0 \Delta p - b_1 \partial_t^{\beta} \Delta p = \tilde{f}\,,
\end{equation*}
\item 
modified Szabo wave equation (fractional Maxwell):
\begin{equation*}
p_{tt}  - a_1\partial_t^{2+\alpha}p - b_0 \Delta p = \tilde{f}\,,
\end{equation*}
\item 
fractional Zener wave equation:
\begin{equation*}
p_{tt}  - a_1\partial_t^{2+\alpha}p - b_0 \Delta p + b_1 \partial_t^{\beta} \Delta p=\tilde{f}\,,
\end{equation*}
\item 
general fractional model:
\[
{\textstyle \sum_{n=0}^N a_n \partial_t^{2+\alpha_n} p - \sum_{m=0}^M b_m \partial_t^{\beta_m} 
\Delta p
=\tilde{f}\,.
}
\]
\end{itemize}
We refer to \cite[Chapter 7]{book_frac} for their well-posednesss analysis.

\subsection{Fractional models of nonlinear acoustics via temperature -- heat flux relation}
Recall that the classical Fourier law \eqref{Fourier} leads to an infinite speed of propagation paradox, which can be amended by substituting it with the Maxwell-Cattaneo law \eqref{MaxwellCattaneo}, allowing for ``thermal waves'' (aka the second sound phenomenon).
However, this can lead to a violation of the 2nd law of thermodynamics and moreover is unable to reproduce the fractional power law frequency dependence of attenuation as relevant in ultrasonics.
Some research in physics has therefore been devoted to an appropriate and physically consistent modeling of the temperature -- heat flux relation.
This can be achieved by ``interpolating'' between \eqref{Fourier} and \eqref{MaxwellCattaneo}, 
using fractional derivatives \cite{CompteMetzler1997,GurtinPipkin1968,Povstenko2011}:
\begin{alignat*}{3}
	\hspace*{-1.5cm}\text{\small (GFE)}\hphantom{III}&&\qquad \qquad  (1+\tau^\alpha \partial_t^\alpha){\boldsymbol{q}}(t) =&&-K \nabla {\vartheta}.\hphantom{{\rt^{1-\alpha}}\partial_t^{1-\alpha} }\\[1mm]
	\hspace*{-1.5cm}\text{\small(GFE I)}\hphantom{II}&& \qquad \qquad(1+\tau^\alpha \partial_t^\alpha)
{\boldsymbol{q}}(t) =&&
\, -K {\rt^{1-\alpha}}\partial_t^{1-\alpha} \nabla {\vartheta}\\[1mm]
	\hspace*{-1.5cm}\text{\small(GFE II)}\, \hphantom{I}&&\qquad \qquad (1+\tau^\alpha \partial_t^\alpha){\boldsymbol{q}}(t) =&&\, -K {\rt^{\alpha-1}}\partial_t^{\alpha-1} \nabla {\vartheta}\\[1mm]
	\hspace*{-1.5cm}\text{\small(GFE III)}\,\, && \qquad \qquad (1+\tau \partial_t){\boldsymbol{q}}(t) =&&\, -K {\rt^{1-\alpha}}\partial_t^{1-\alpha} \nabla {\vartheta}. 
\end{alignat*}
Using these in place of \eqref{MaxwellCattaneo} leads to the models
\begin{equation}\label{fracJMGTeqs}\begin{aligned}
&\tau^\alpha \partial_t^{\alpha} \psi_{tt}+(1-\kappa\psi_t)\psi_{tt}-c^2 \Delta \psi -\tau^\alpha c^2 \partial_t^{\alpha} \Delta \psi- \delta \Delta\psi_{t}+ \partial_t |\nabla \psi|^2=0\\
&\tau^\alpha \partial_t^{\alpha} \psi_{tt}+(1-\kappa\psi_t)\psi_{tt}-c^2 \Delta \psi -\tau^\alpha c^2 \partial_t^{\alpha} \Delta \psi- \delta \partial_t^{2-\alpha} \Delta\psi+ \partial_t |\nabla \psi|^2=0\\
&\tau^\alpha \partial_t^{\alpha} \psi_{tt}+(1-\kappa\psi_t)\psi_{tt}-c^2 \Delta \psi -\tau^\alpha c^2 \partial_t^{\alpha} \Delta \psi- \delta \partial_t^{\alpha} \Delta\psi+\partial_t |\nabla \psi|^2=0\\
&\tau \psi_{ttt}+(1-\kappa\psi_t)\psi_{tt}-c^2 \Delta \psi -\tau c^2 \Delta \psi_{t}- \delta \partial_t^{2-\alpha}\Delta\psi+ \partial_t |\nabla \psi|^2=0
\end{aligned}\end{equation}
that have been analyzed in \cite{fracJMGT}.

\subsection{Fractional derivatives}\label{sec:fracderiv}
To clarify what we mean by $\partial_t^\alpha$ in the models above, we provide a brief intermission on some of the most commonly used fractional derivative concepts.

With the \Blue{Abel fractional integral operator}
\[ 
I_{t_0}^\gamma f(x) = \frac{1}{\Gamma(\gamma)}\int_{t_0}^t \frac{f(s)}{(t-s)^{1-\gamma}}\,ds,
\]
a fractional (time) derivative can be defined by either
\[
\begin{aligned}
&{}^{RL}_{t_0} D^\alpha_t f  = \frac{d\ }{dt} I_{t_0}^{1-\alpha} f\quad \mbox{\Blue{Riemann-Liouville derivative}}
\text{ or }\\
&{}^{DC}_{t_0} D^\alpha_t f = I_{t_0}^{1-\alpha} \frac{df}{ds} \quad \mbox{\Blue{Djrbashian-Caputo derivative}}
\end{aligned}
\]
These are nonlocal operators and have a definite starting point $t_0$. 
The R-L derivative is defined on a larger function space that the D-C one. However, the R-L derivative of a constant is nonzero; it even exhibits a singularity at initial time $t_0$.
This is why the D-C derivative, which maps constants to zero, is often preferred in physical applications.
Some recent books on fractional PDEs (also containing many classical references on fractional differentiation) are \cite{Jin:2021,book_frac,KubicaRyszewskaYamamoto:2020}.

The nonlocal and causal character of these derivatives provides the models in which they are used with a ``memory''.
In the context of inverse problems, it is important to note that therefore initial values are tied to later values and hence can be better reconstructed backwards in time than in integer derivative models.

Some of the challenges in the analysis of fractional PDEs, in particular in deriving energy estimates, are as follows.
\begin{itemize}
\item The \Red{chain rule} identity $u \, u_t = \frac12 \frac{d}{dt} u^2$ is lost;
\item Likewise, there is no \Red{product rule} of differentiation.
\end{itemize}
{``Substitutes'' for these devices in the fractional derivative case are}
\begin{itemize}
\item the chain rule inequality
\begin{equation*}
	\Blue{{w(t)}\textup{D}_t^{\alpha}w(t)\geq \tfrac12(\textup{D}_t^{\alpha} w^2)(t)}
	\end{equation*}
for $w\in W^{1,1}(0,T)$ \cite{Alikhanov:11};
\item coercivity of the Abel integral operator
	\begin{equation*}
	\Blue{\int_0^T \langle \textup{I}^{1-\alpha} w(s),  w(s) \rangle \,ds \geq \cos ( \tfrac{\pi(1-\alpha)}{2} ) \| w \|_{H^{-(1-\alpha)/2}(0,T)}^2} 
	\end{equation*}	 
for $w\in H^{-(1-\alpha)/2}(0,T)$ \cite{Eggermont1987,VoegeliNedaiaslSauter2016};
\item the Kato-Ponce inequality 
\begin{equation*}
\Blue{\|f\,g\|_{W^{\rho,r}(0,T)}\lesssim 
\|f \|_{W^{\rho,p_1}(0,T)} \|g\|_{L^{q_1}(0,T)}
+ \| f \|_{L^{p_2}(0,T)} \|g\|_{W^{\rho,q_2}(0,T)}
}
\end{equation*}
for $f,\,g$ such that the right hand side is finite, $0\leq\rho\leq\overline{\rho}<1$, 
$1<r<\infty$, $p_1$, $p_2$, $q_1$, $q_2\in(1,\infty]$, with $\frac{1}{r}=\frac{1}{p_i}+\frac{1}{q_i}$, $i=1,2$; see, e.g., \cite{GrafakosOh2014}.
\end{itemize}

\section{Describing nonlinear wave propagation in frequency domain by 
multiharmonic expansions}

Motivated by the characteristic appearance of higher harmonics in nonlinear acoustics, cf. Figure~\ref{fig:harmonics}, we now provide a mathematical formulation of this physical fact.

To this end, we recall the linear wave equation 
\begin{equation}\label{lin_wave}
u_{tt}- c^2\Delta u  = r,
\end{equation}
where using a harmonic excitation $r(x,t)=\Re(\hat{r}(x)e^{\imath \omega t})$ 
and a harmonic ansatz for $u$ 
\[
u(x,t)= \Re\left(\hat{u}(x) e^{\imath \omega t}\right)
\]
leads to the Helmholtz equation
\[
-\omega^2 \hat{u}- c^2\Delta \hat{u}  = \hat{r}.
\]
Replacing \eqref{lin_wave} by the Westervelt equation
\begin{equation}\label{Westervelt_1}
u_{tt}- c^2\Delta u - b \Delta u_t =
\kappa (u^2)_{tt} + r
\end{equation}
and still using a harmonic excitation $r(x,t)=\Re(\hat{r}(x)e^{\imath \omega t})$, we expect a response at multiples of the fundamental frequency $\omega$. 
We therefore use a \Blue{multiharmonic expansion} of $u$ 
\begin{equation}\label{ansatz_multiharmonic}
u(x,t)
= \Re\left(\sum_{k=1}^\infty \hat{u}_k(x) e^{\imath k \omega t}\right),
\end{equation}
which is mathematically justified by completeness of the system $(e^{\imath k \omega \cdot})_{k\in\mathbb{N}}$ in 
$L^2(0,T)$
as well as the fact that indeed a periodic solution to the Westervelt equation exists \cite{periodicWestervelt}.

\begin{theorem}\label{th:wellposedness}
	For $b,c^2,\beta,\gamma,T>0$, $\kappa\in L^\infty(\Omega)$, there exists $\rho>0$ such that for all $r\in L^2(0,T;L^2(\Omega)$ with $\|r\|_{L^2(0,T;L^2(\Omega))}\leq\rho$ there exists a unique solution 
\[
	\begin{aligned}
	 u \in \, X:=H^2(0,T;L^2(\Omega))\cap H^1(0;T;H^{3/2}(\Omega)) \cap L^2(0;T;H^2(\Omega))
	\end{aligned}
\]
 of
	\begin{equation*} 
		\begin{aligned}
		\begin{cases}
		u_{tt}-c^2\Delta u - b\Delta u_t = \kappa(x) (u^2)_{tt} +r \quad \mbox{ in }\Omega\times(0,T), \\[2mm]
		\beta u_t+\gamma u+\partial_\nu u=0 \quad \mbox{ on } \partial\Omega\times(0,T),\\[2mm]
		u(0)=u(T)\,, \ u_t(0)=u_t(T) \quad \mbox{ in }\Omega,
		\end{cases}
		\end{aligned}
		\end{equation*}
and the solution fulfills the estimate
	\begin{equation*}
	\begin{aligned}
\|u\|_X\leq \tilde{C} \|r\|_{L^2(0,T;L^2(\Omega))}.
	\end{aligned}
	\end{equation*}
\end{theorem}

Inserting the ansatz \eqref{ansatz_multiharmonic} into \eqref{Westervelt_1} leads to the coupled system
\begin{equation}\label{multiharmonics}
\hspace*{-0.7cm}\begin{aligned}
&m=1:&& -\omega^2 \hat{u}_1-(c^2+\imath\omega b) \Delta \hat{u}_1 = \hat{r}
\textcolor{grey}{ \ -\frac{\kappa}{2}\omega^2
\sum_{k=3:2}^\infty\overline{\hat{u}_{\frac{k-1}{2}}} \hat{u}_{\frac{k+1}{2}}}
\\
&m=2,3\ldots:&& -\omega^2 m^2 \hat{u}_m-(c^2+\imath\omega m b) \Delta \hat{u}_m 
= -\frac{\kappa}{4}\omega^2 m^2 \sum_{\ell=1}^{m-1} \hat{u}_\ell \hat{u}_{m-\ell} 
 \\&&&\hspace*{5cm}\textcolor{grey}{ 
-\frac{\kappa}{2}\omega^2 m^2 \sum_{k=m+2:2}^\infty\overline{\hat{u}_{\frac{k-m}{2}}} \hat{u}_{\frac{k+m}{2}}}.
\end{aligned}
\end{equation}
Note that the quadratic nonlinearity corresponds to a discrete autoconvolution in time.
Skipping the grey terms corresponds to the simplification of dropping $\Re$ in \eqref{ansatz_multiharmonic} and leads to a triangular system that can be solved by forward substitution in a sequence of Helmholtz problems.

\section{Inverse problems}\label{sec:inverse}
In this section we mainly focus on an innovative quantitative imaging methodology, based on some of the models in the previous sections, 
namely acoustic nonlinearity parameter tomography ANT.
The key motivation for this comes from the experimentally proven qualification of the B/A parameter in \eqref{Kuznetsov}, \eqref{Westervelt}, \eqref{JMGT} or \eqref{BCBJ} as a characterizing and distinguishing material property of biological tissues. 
Thus, viewing $\kappa=\tfrac{1}{\varrho c^2}(\tfrac{B}{2A}+1)$ as a spatially varying coefficient in the Westervelt equation, it can be used for medical imaging.
While the physical and engineering literature on ANT has already been developping over many years, with initializing work going back to the 1980's and 90's, \cite{Bjorno1986,BurovGurinovichRudenkoTagunov1994,Cain1986,IchidaSatoLinzer1983,ZhangChenYe1996},
mathematical results, as they are crucial for enabling any kind of advanced medical imaging technology, have here only been started off recently. 
In particular, we point to \cite{AcostaUhlmannZhai2021,UhlmannZhang2023} for results on uniqueness from the  Neumann-Dirichlet map in the Westervelt equation and to \cite{Eptaminitakis:Stefanov:2023}, which investigates and highlights the usefulness of geometric optics solutions to Westervelt equation for this imaging problem. 
In our own work \cite{nonlinearity_imaging_Westervelt,nonlinearity_imaging_fracWest,nonlinearity_imaging_JMGT,nonlinearity_imaging_2d,nonlinearity_imaging_both}
we mainly focused on uniqueness from a single boundary observation and on convergence of Newton's method in the Westervelt equation as well as
in \cite{nonlinearityimaging} on linearized uniqueness and conditional stability of the inverse problem for the BCBJ model \eqref{BCBJ}.

With the Westervelt equation as a forward model, the inverse problem of nonlinearity parameter imaging reads as follows.\\[1ex]
\textit{Reconstruct \Red{$\kappa(x)$} in 
\begin{equation}\label{eqn:Westervelt_init_D_intro}
\begin{aligned}
&\bigl(u-\Red{\kappa(x)}u^2\bigr)_{tt}-c_0^2\Delta u + D u = r \quad
\mbox{ in }\Omega\times(0,T)\\
\partial_\nu u+\gamma u&=0 \mbox{ on }\partial\Omega\times(0,T),\quad
u(0)=0, \quad u_t(0)=0 \quad \mbox{ in }\;\Omega
\end{aligned}
\end{equation}
(with excitation $r$) from observations 
\begin{equation*}
g= u \quad \mbox{ on }\Sigma\times(0,T)
\end{equation*}
on a transducer array $\Sigma\subset \overline\Omega$, which we assume to be a surface or  collection of discrete points.
}

The attenuation term $D$ can be defined by one of the following common fractional models.
\[\begin{aligned}
&\text{Caputo-Wismer-Kelvin: }
\qquad D = -b \Delta \partial_t^\beta  
\quad\mbox{ with }\beta\in[0,1], \ \ b\geq0 
\\
&\text{fractional Zener: }
\qquad D = \aone  \partial_t^{2+\alpha} - \bone  \Delta \partial_t^{\beta}   
\quad \mbox{ with }\aone >0, \ \bone \geq \aone c^2, \  1\geq \beta\geq\alpha>0,
\\
&\text{space fractional Chen-Holm: } 
\qquad D =  b (-\Delta)^{\betil} \partial_t  
\quad\mbox{ with }\betil\in[0,1], \ b\geq0. 
\end{aligned}\]

Since also the speed of sound and the attenuation coefficient depend on the underlying tissue type, they can be considered as additional unknowns in a multicoefficient inverse problem, see, e.g., \cite{nonlinearity_imaging_both,nonlinearity_imaging_JMGT}.

{Chances and challenges} of this inverse problem result from 
\begin{itemize}
\item the difficulties in the analysis of the forward model (as well as its numerics) that have already been pointed to in Section~\ref{sec:models}; 
\item the fact that the unknown coefficient $\kappa(x)$ appears in the nonlinear term;
\item $\kappa(x)$ being spatially varying whereas the data $g(t)$ is in the
``orthogonal'' time direction; 
this is well known to lead to severe 
ill-conditioning of the inverse problem.
\end{itemize}

On the other hand, nonlinearity helps by adding information.
In particular the generation of higher harmonics and the fact that $\kappa$ appears in each of the equations in \eqref{multiharmonics} shows that information multiplies due to nonlinearity.
As a consequence, it is in fact possible to achieve enhanced uniqueness results for the inverse problem formulated above, which we here illustrate by a uniqueness result for a general spatially variable coefficient from a single observation, while in the linear setting, this type of data would at most provide identifiability of piecewise constant coefficients. 

\subsection{Linearized and local nonlinear uniqueness}

To this end, recalling \eqref{multiharmonics} we reformulate 
the inverse problem in frequency domain, using the abbreviations 
\[
\vec{u}=(\hat{u}_j)_{j\in\mathbb{N}} \quad
B_m(\vec{u})= \frac14 \sum_{\ell=1}^{m-1} \hat{u}_\ell \hat{u}_{m-\ell} 
+\frac12 \sum_{k=m+2:2}^\infty\overline{\hat{u}_{\frac{k-m}{2}}} \hat{u}_{\frac{k+m}{2}},
\] 
defining the forward operator by 
\[
F_m(\Red{\kappa},\Blue{\vec{u}}):=
\left(\begin{array}{l}
-\bigl(\omega^2 m^2 +(c^2+\imath\omega m b) \Delta\bigr) \Blue{\hat{u}_m}
+\omega^2 m^2 \, \Red{\kappa}\, B_m(\Blue{\vec{u}})\\
\text{tr}_\Sigma \Blue{\hat{u}_m}
\end{array}\right)
\]
and the data by
\[
y_m:= \left(\begin{array}{l}
\hat{r} \text{ if }m=1 \ / \ 0 \text{ if } m\geq2\\
\hat{g}_m
\end{array}\right).
\]
This allows us to write the inverse problem as a nonlinear operator equation
\[
\vec{F}(\Red{\kappa},\Blue{\vec{u}}) = y.
\]
The linearized inverse problem is defined by considering the 
linearization of the forward operator
\[
\vec{F}'(\kappa,\vec{u})(\Red{\uld{\kappa}},\Blue{\uld{\vec{u}}}) \approx F(\kappa,\vec{u})-y,
\]
where 
\[
\begin{aligned}
&F_m'(\kappa,\vec{u})(\Red{\uld{\kappa}},\Blue{\uld{\vec{u}}})=\\
&
\left(\begin{array}{l}
-\bigl(\omega^2 m^2 +(c^2+\imath\omega m b) \Delta\bigr) \Blue{\uld{\hat{u}}_m}
+\omega^2 m^2 \, \kappa\, B_m'(\vec{u})\Blue{\uld{\vec{u}}}
+\omega^2 m^2 \, \Red{\uld{\kappa}}\, B_m(\vec{u})
\\
\text{tr}_\Sigma \Blue{\uld{\hat{u}}_m}
\end{array}\right)
\end{aligned}
\]
and in particular
\[
\begin{aligned}
&F_m'(0,\vec{u})(\Red{\uld{\kappa}},\Blue{\uld{\vec{u}}})=
\left(\begin{array}{l}
-\bigl(\omega^2 m^2 +(c^2+\imath\omega m b) \Delta\bigr) \Blue{\uld{\hat{u}}_m}
+\omega^2 m^2 \, \Red{\uld{\kappa}}\, B_m(\vec{u})
\\
\text{tr}_\Sigma \Blue{\uld{\hat{u}}_m}
\end{array}\right)
\end{aligned}
\]

This allows to formulate the following linearized uniqueness result, that has been obtained for the Westervelt equation in \cite{nonlinearity_imaging_2d} and for the JMGT equation in \cite{nonlinearity_imaging_JMGT_freq}.

\begin{theorem}
The homogeneous linearised 
(at $\kappa=0$, $\vec{u}=\phi(x)\vec{\psi}$ with $B_m(\vec{\psi})\not=0$, $m\in\mathbb{N}$)\\
inverse problem of nonlinearity coefficient imaging in frequency domain only has the trivial solution, that is, $\vec{F}'(0,\vec{u})$ is injective.
\end{theorem}

In the JMGT model, this even enables a (local) nonlinear uniqueness result
\cite{nonlinearity_imaging_JMGT_freq}.
\begin{theorem}
For $(\kappa,\vec{u})$, $(\tilde{\kappa},\tilde{\vec{u}})$ in a sufficiently small 
$H^s(\Omega)\times h^1(L^2(\Omega))$ neighborhood of $(0, \phi(x)\vec{\psi})$ with $s\in(\frac12,1]$ 
the uniqueness result 
\[
\vec{F}(\kappa,\vec{u})=\vec{F}(\tilde{\kappa},\tilde{\vec{u}}) \quad \Rightarrow \quad 
\kappa=\tilde{\kappa}\text{ and } \vec{u}=\tilde{\vec{u}}
\]
holds.
\end{theorem}

{}
\subsection{Some reconstruction results from two harmonics}
\begin{figure}
\includegraphics[width=0.19\textwidth]{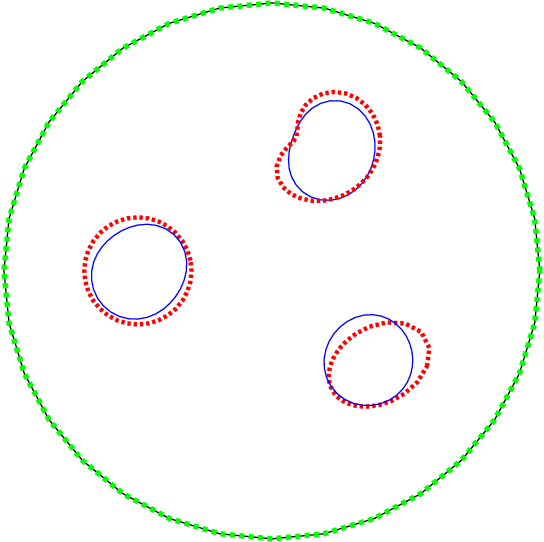}
\includegraphics[width=0.19\textwidth]{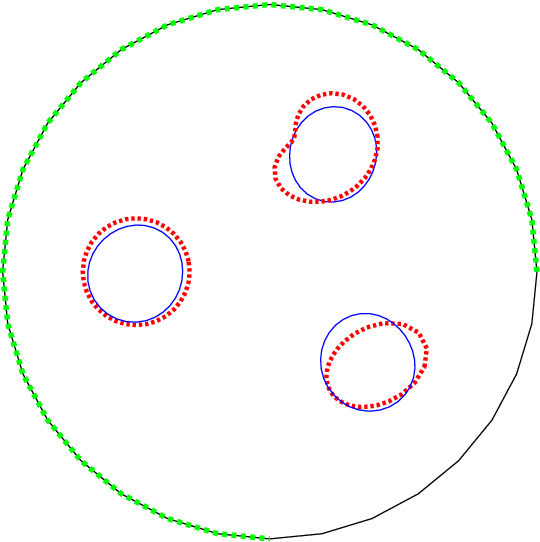}
\includegraphics[width=0.19\textwidth]{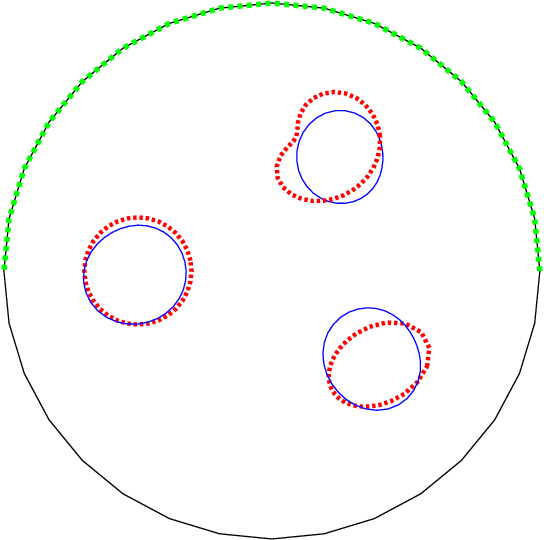}
\includegraphics[width=0.19\textwidth]{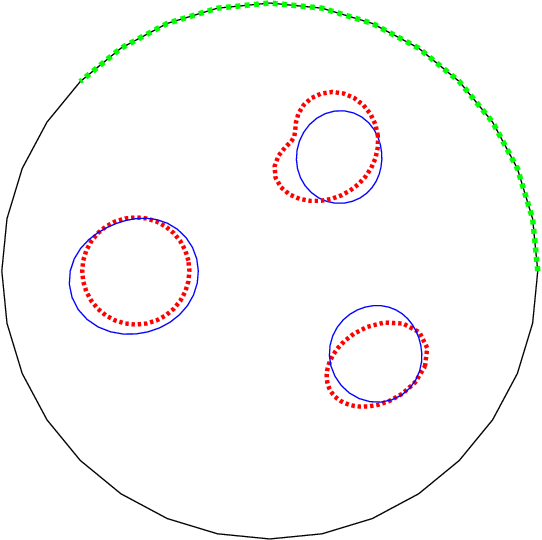}
\includegraphics[width=0.19\textwidth]{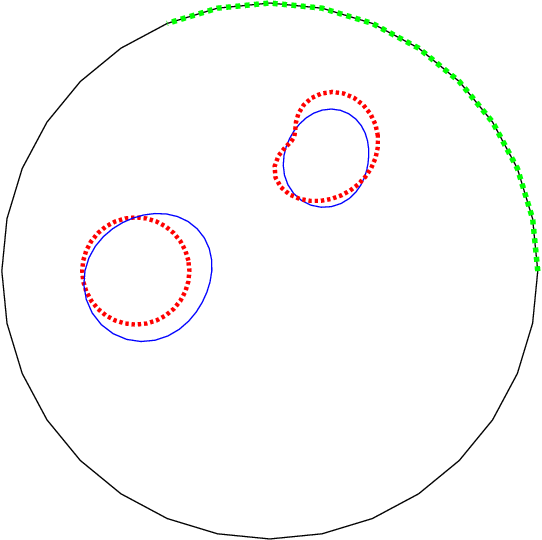}
\\[2ex]
(a) $\frac{\alpha}{2\pi}=1$ \hspace*{0.015\textwidth}
(b) $\frac{\alpha}{2\pi}=0.75$ \hspace*{0.015\textwidth}
(c) $\frac{\alpha}{2\pi}=0.5$ \hspace*{0.015\textwidth}
(d) $\frac{\alpha}{2\pi}=0.4$ \hspace*{0.015\textwidth}
(e) $\frac{\alpha}{2\pi}=0.3$ \hspace*{-0.01\textwidth}
\caption{reconstructions from limited data: measurements are taken on the green part of the boundary;  $\alpha$ corresponds to the observation angle.
\label{fig:3objects}}
\end{figure}

\begin{figure}
\includegraphics[width=0.19\textwidth]{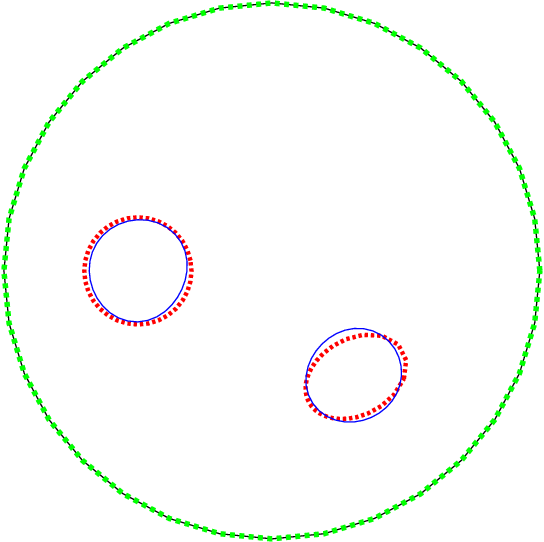}
\hspace*{0.02\textwidth}
\includegraphics[width=0.19\textwidth]{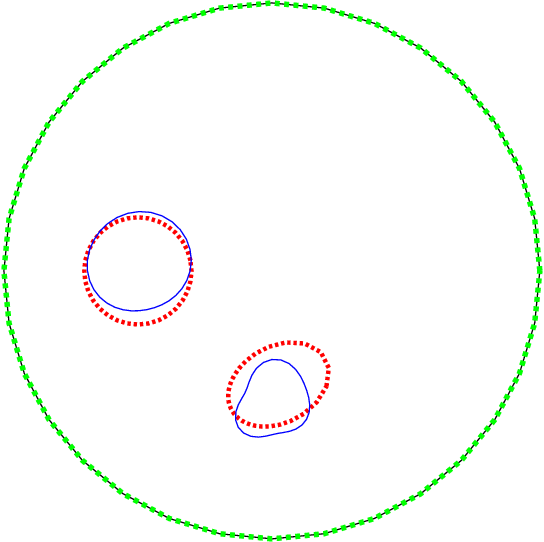}
\hspace*{0.02\textwidth}
\includegraphics[width=0.19\textwidth]{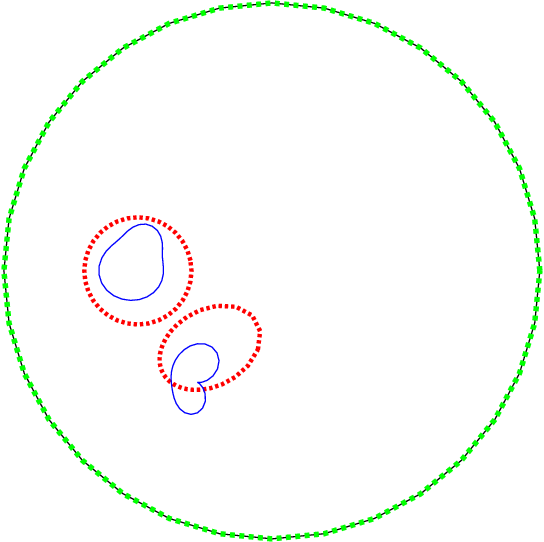}
\hspace*{0.02\textwidth}
\includegraphics[width=0.19\textwidth]{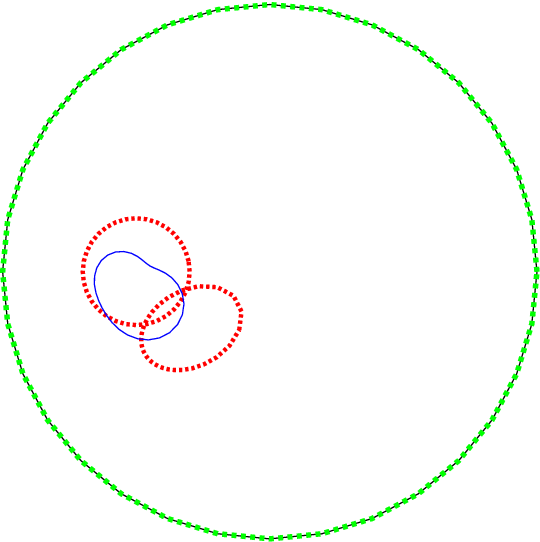}
\\[2ex]
(a) $\frac{\theta}{2\pi}=0.3$ \hspace*{0.05\textwidth}
(b) $\frac{\theta}{2\pi}=0.2$ \hspace*{0.05\textwidth}
(c) $\frac{\theta}{2\pi}=0.1$ \hspace*{0.05\textwidth}
(d) $\frac{\theta}{2\pi}=0.09$ \hspace*{-0.01\textwidth}
\caption{{Reconstruction of two inclusions at different distances from each other}; $\theta$ corresponds to the distance in polar coordinates.
\label{fig:2objects}}
\end{figure}

\begin{figure}
\includegraphics[width=0.19\textwidth]{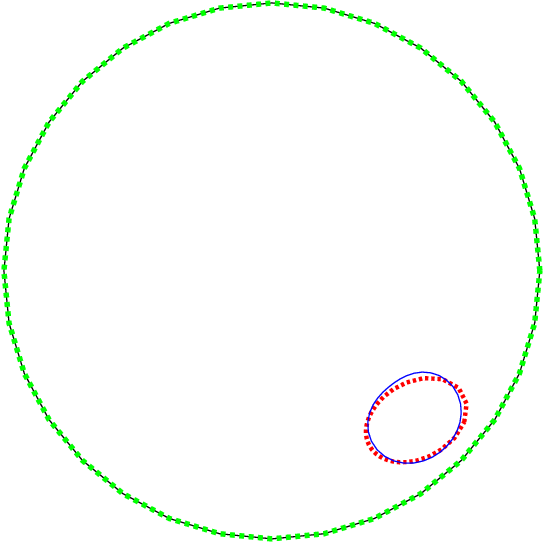}
\hspace*{0.02\textwidth}
\includegraphics[width=0.19\textwidth]{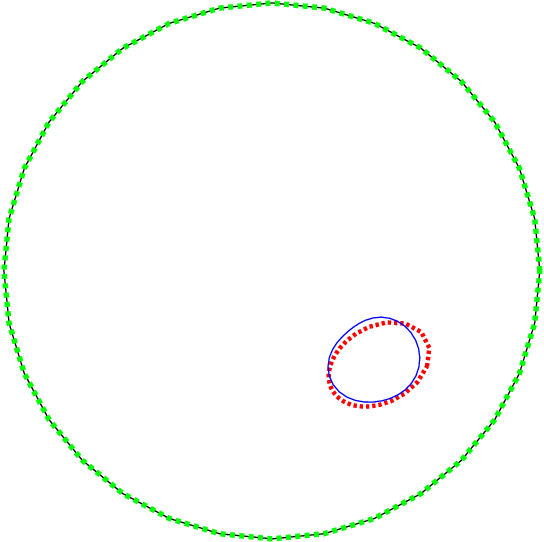}
\hspace*{0.02\textwidth}
\includegraphics[width=0.19\textwidth]{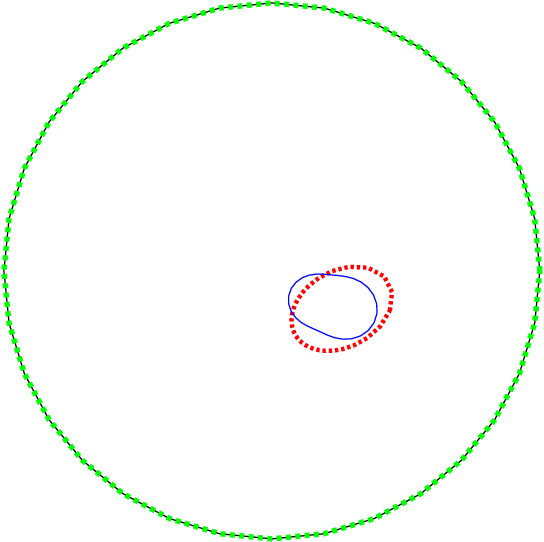}
\caption{{Reconstruction of one inclusion at different distances from the boundary}
\label{fig:1object}}
\end{figure}

Finally we provide some images under the (practically often relevant) assumption of a piecewise constant nonlinearity coefficient, using only the information provided by the first two lines  $m\in\{1,2\}$ of the multiharmonic system \eqref{multiharmonics}.
The reconstruction was carried out by parametrizing the boundary of the interface curves and applying Newton's method to the nonlinear system that arises from matching measurements to simulations. Starting curves for Newton's method were found by applying a sparsity enhancing point source insertion algorithm \cite{BrediesPikkarainen13} and extending the discovered point sources to disks whose radii are determined by a mean value property of solutions to the Helmholtz equation.
While the geometry in Figures~\ref{fig:3objects}, \ref{fig:2objects}, \ref{fig:1object} is clearly idealized, we mimic a realistic measurement scenario by adding 1\% random noise to the simulated data.
For details we refer to \cite{nonlinearity_imaging_2d} where also the pictures shown here are taken from. 

In Figure~\ref{fig:3objects} we consider three objects and study the impact of limited data on the reconstruction quality, which clearly decreases under loss of information. Remarkably, still all three objects are found with less than half of the full view and only upon reduction to less than a third of the data one of the three objects gets invisible.

Figure~\ref{fig:2objects} illustrates a study of the influence of the distance bewteen two objects on their reconstructability. Obviously imaging both of them becomes harder as they move closer to each other. Still, up to the point where the two inclusions touch they can be identified as two objects.

Finally, in Figure~\ref{fig:1object} we visualize how attenuation compromises the reconstruction as the object moves further away from the boundary on which observations are taken.


\begin{ack}
The author thanks Benjamin Rainer, University of Klagenfurt and Austrian Institute of Technology, for providing the illustrative images in Figure~\ref{fig:harmonics}.
\end{ack}

\begin{funding}
This research was funded in part by the Austrian Science Fund (FWF) 
[10.55776/P36318]. 
\end{funding}










\end{document}